\documentclass[a4paper]{article}

\usepackage[table]{xcolor}
\usepackage{colortbl}
\usepackage{geometry}

\usepackage[nomarkers, nolists,figuresonly]{endfloat}
\usepackage{booktabs}
\usepackage{longtable}
\usepackage{array}
\usepackage{natbib}
\usepackage{amsmath,amssymb}
\usepackage{multirow}
\usepackage{wrapfig}
\usepackage{float}
\usepackage{pdflscape}
\usepackage{tabu}
\usepackage{threeparttable}
\usepackage[normalem]{ulem}

\usepackage{authblk}

\usepackage{setspace}
\spacing{1.5}
\usepackage{amsthm}
\newtheorem{theorem}{Theorem}[section]
\newtheorem{lemma}{Lemma}[section]

\title{Moment-based Estimation of Mixtures of Regression Models}

\author[1]{Claus Thorn Ekstr{\o}m}
\author[1]{Christian Bressen Pipper}
\affil[1]{Section of Biostatistics, Department of Public Health,\\ University of Copenhagen, \O ster Farimagsgade 5 B, DK-1014 Copenhagen K, Denmark}

\begin{document}


\maketitle


\begin{abstract} {Finite mixtures of regression models provide a
    flexible modeling framework for many phenomena. Using moment-based
    estimation of the regression parameters, we develop unbiased
    estimators with a minimum of assumptions on the mixture
    components.  In particular, only the average regression model for
    one of the components in the mixture model is needed and no
    requirements on the distributions. The consistency and asymptotic
    distribution of the estimators is derived and the proposed method
    is validated through a series of simulation studies and is shown
    to be highly accurate. We illustrate the use of the moment-based
    mixture of regression models with an application to wine quality
    data.}
  {mixture regression model; moment estimation; estimating
    equations; robust inference}
\end{abstract}

\section{Introduction}
\label{sec1}

The class of finite mixtures of regression models provides a flexible
approach to model a wide range of phenomena and to handle non-standard
data analysis problems such as excess zeros or heterogeneity
\citep{mcla:peel:2000}. Mixture models can also be used to accommodate data
that is "contaminated" due to poor quality data, laboratory errors, or
situations where data originates from multiple sources.  In
particular, zero-inflated regression models, and hurdle models can be
considered special cases of the class of finite mixture of regression
models with two components and both of these types of models see
frequent use \citep{lamb:1992, kwag:appr:2016}.

Finite mixture regression models are relevant for analyzing many
problems such as sudden-infant-death-syndrome \citep{dalr:etal:2003},
HIV-risk reduction trials \citep{hu:etal:2011}, medical care
\citep{deb:triv:1997}. However, problems relevant for finite mixtures are
especially common in genomics where for example population admixture
(certain subgroups of the population do not segregate a
phenotype-influencing mutation), or in the presence of
gene-environment interactions, where the effect of some genes are
never triggered because the individual is living under specific
environmental conditions.

Other recent approaches of using finite mixture models in genetics
include the paper by \citet{xu:etal:2015} where they compare different model
types (standard parametric, non-parametric, zero-inflated and hurdle
models) to microbiome data in order to assess and infer the best model
for these types of data. They find that mixture models in general fit
the data best but that the choice of parametric model (Poisson or
negative-binomial) can have substantial impact on the results and on
the convergence on the estimation algorithm.

Typically, parameters in mixture regression models are estimated by
specifying a fully parametric model for each of the components. This
leads to efficient estimates provided that the parametric models are
specified correctly.  In this article we will consider the class of
finite mixture regressions where the predictors are influencing
exactly one of the regression model components. Within this framework
we wish to estimate the association between the predictors and the
outcome. We place no restriction on the distribution of the components
but only require that the relationship between the predictors and the
outcome can be modeled as a linear model. Moment-based estimators
allow us to obtain unbiased estimates of the regression parameters
with a minimum of assumptions. In particular, we do not require a
specification of the distribution, and we only need to specify the
average regression model for one of the components in the mixture
model.

The recent paper by \citet{kong:etal:2017} uses an approach for
identifying interaction that most closely resembles the approach in
this manuscript. They also model the first and second
moments but they are focused on interaction in particular and do not
use a moment-based method for estimation.

The manuscript is organized as follows: The next section explains the
problem, derives the moment-based estimators, and proves the
large-sample properties relevant for inference. In section 3, we
present properties of the proposed method for moment-based mixture of
regression models and use simulations to compare the proposed
moment-based mixture regression model to the analogous Gaussian
mixture model.  Finally, we apply the moment-based mixture of
regression models to a dataset of wine quality before we discuss the
findings along with possible extensions. The approaches presented in
this manuscript are available in the R package \texttt{mommix} which
can be found on github at \texttt{www.github.com/ekstroem/mommix}.

\section{Methods}\label{methods}

Let \(Y\) be given as a mixture of two distributions, that is;

\[Y = PY_1 + (1-P)Y_2,\] where \(\text{Pr}(P=1)=1-\text{Pr}(P=0)=p\) for
some \(p\in(0,1)\) and \(Y_1\) and \(Y_2\) are independent stochastic
variables from the two underlying distributions. We are interested in
modeling the relationship between a set of \(m\) predictors represented
by the design matrix \(X\) and \(Y_1\) and can think of \(Y_2\) as a
contamination of the response.

For the distribution of interest we assume that
\(E(Y_1|X) =\mu_1+ \beta^{T}X\), and \(V(Y_1|X) = \sigma_1^2\), but for
the other component we only require that \(Y_2\perp X\) and that it has
well-defined mean and variance \(E(Y_2) = \mu_2\), and
\(V(Y_2) = \sigma_2^2\).

The mean of \(Y\) is

\begin{equation} \label{firstorder}
E(Y|X) =p\cdot(\mu_{1}+\beta^{T}X) + (1-p)\cdot\mu_2,
\end{equation}

where \(\beta\in \Bbb R^m\) is the parameter vector of interest. Simple
algebra yields that the second order moment is

\begin{equation}\label{secorder}
E(Y^2|X) = p\cdot( (\mu_{1}+\beta^{T}X)^{2} + \sigma_1^2) + (1-p)\cdot(\mu_{2}^{2} + \sigma_2^2).
\end{equation}

Note that the formula for the second order moment has a form
corresponding to a multiple linear regression model with a mean that can
be reparameterized as
\(2\cdot\mu_{1}\cdot p\cdot(\beta^{T} X)+p\cdot(\beta^{T} X)^2 + \tilde\alpha\).
Thus, if \(X\) is univariate we can obtain estimates of
\(2\cdot\mu_{1}\cdot p\cdot \beta\) and \(p\cdot\beta^2\) directly from
least-squares regression of the second order moment and estimate
\(p\cdot\beta\) from least-squares regression of the first moment. These
estimates can be combined to extract individual estimates of
\(\mu_{1}\), \(\beta\), and \(p\). In the following we formalize and
extend this idea.

Let \(\lambda_1=p\cdot\beta\) and set \(\eta_{i}=\lambda_{1}^{T}X_{i}\)
(the linear predictor), \(\lambda_{2}=2\cdot\mu_{1}\), and
\(\lambda_{3}=p^{-1}\). Then equation (\ref{secorder}) may be rewritten
as

\begin{equation*}
E(Y_{i}^2|X_{i})=\tilde{\alpha}+\lambda_{2}\cdot \eta_{i}+\lambda_{3}\cdot \eta_{i}^{2},
\end{equation*}

where the intercept \(\tilde{\alpha}\) contains contributions from both
\(Y_1\) and \(Y_2\) but holds no direct information about \(X\).

From this parameterization we have that
\(\beta=\lambda_1\cdot\lambda_3\), \(\mu_1 = \lambda_2/2\), and
\(p=\frac1{\lambda_3}\). Thus, if we can estimate \(\lambda_1\),
\(\lambda_2\), and \(\lambda_3\) then we can also estimate the mixture
proportion and all mean parameters related to the distribution of
interest.

\subsection{Inference and large sample
properties}\label{inference-and-large-sample-properties}

Note that \(\lambda_1=p\cdot\beta\) corresponds to the regression
parameter vector when we regress \(Y\) on \(X\). Consequently by least
squares regression we can obtain a closed form estimator
\(\hat{\lambda}_{1}\) of \(\lambda_{1}\). Specifically with
\(X_i\in \mathbb{R}^{m}\) and \(Y_i\in \mathbb{R}\) denoting the data
from the \(i\)th individual and with
\(\tilde{X}=(\tilde{X}_1,\ldots,\tilde{X}_n)^{T}\),
\(\tilde{X}_{i}=(1,X_{i}^{T})^{T}\), and \(Y=(Y_1,\ldots,Y_n)\) we have:

\begin{equation}
\hat{\lambda}_{1}=A(m)(\tilde{X}^{T}\tilde{X})^{-1}\tilde{X}^{T}Y,
\end{equation}

where \(A(m)\) is a \(m\times(m+1)\) matrix projecting an \(m+1\) vector
to its last \(m\) components.

The residual \(Y_{i}-E(Y_{i}|\tilde{X}_{i})\) is not ensured to follow a
Gaussian distribution --- even if the original distributions of \(Y_1\)
and \(Y_2\) were Gaussian --- but we may still resort to a
characterization of \(\sqrt{n}(\hat{\lambda}_{1}-\lambda_{1})\) as a sum
of i.i.d. zero mean terms ensuring consistency and asymptotic normality
by means of the central limit theorem. Specifically the characterization
is as follows:

\begin{equation}\label{iid1}
\sqrt{n}(\hat{\lambda}_{1}-\lambda_{1})=A(m)(\frac{1}{n}\tilde{X}^{T}\tilde{X})^{-1}\frac{1}{\sqrt{n}}\sum_{i=1}^{n}\tilde{X}_{i}\{Y_{i}-E(Y_{i}|\tilde{X}_{i})\}=\frac{1}{\sqrt{n}}\sum_{i=1}^{n}\varepsilon_{i}+o_{P}(1)
\end{equation}

with

\begin{equation*}
\varepsilon_{i}=A(m)E(\tilde{X}_{i}\tilde{X}_{i}^{T})^{-1}\tilde{X}_{i}\{Y_{i}-E(Y_{i}|\tilde{X}_{i})\}.
\end{equation*}

Recall that equation (\ref{secorder}) may be rewritten as

\begin{equation*}
E(Y_{i}^2|X_{i})=\tilde{\alpha}+\lambda_{2}\cdot \eta_{i}+\lambda_{3}\cdot \eta_{i}^{2}
\end{equation*}

with \(\eta_{i}=\lambda_{1}^{T}X_{i}\) (the linear predictor),
\(\lambda_{2}=2\cdot\mu_{1}\), and \(\lambda_{3}=p^{-1}\). From this
observation and the fact that we have an estimator of \(\lambda_{1}\) it
seems only natural to estimate \(\lambda_{2}\) and \(\lambda_{3}\) by a
multiple linear regression of \(Y_i^2\) on
\(\hat{\eta}_{i}=\hat{\lambda}_{1}^{T}X_{i}\). However, since the
variance of \(Y_{i}^{2}\) given \(X_{i}\) depends on \(\eta_{i}\) and
\(\eta_{i}^{2}\) it is more appropriate to estimate \(\lambda_{2}\) and
\(\lambda_{3}\) using a weighted regression, where the weight is a
function of \(\eta_{i}\).

Specifically let \(w\) be a weight function, put \(w_{i}=w(\eta_{i})\),
\(\hat{w}_{i}=w(\hat{\eta}_{i})\), and let \(\mathbb{P}_{n}^{w}\) denote
the weighted empirical mean, that is for \(f_{i}=f(Y_{i},X_{i})\) \[
\mathbb{P}_{n}^{w}(f)=\frac{\sum_{i=1}^{n}f_{i}w_{i}}{\sum_{i=1}^{n}w_{i}}.
\] Also define

\begin{align*}
&a_{1}^{n}(\lambda_{1})=\mathbb{P}_{n}^{w}(\eta^{4})-\mathbb{P}_{n}^{w}(\eta^{2})^{2},\\
&a_{2}^{n}(\lambda_{1})=\mathbb{P}_{n}^{w}(\eta^{3})-\mathbb{P}_{n}^{w}(\eta)\mathbb{P}_{n}^{w}(\eta^{2}),\\
&a_{3}^{n}(\lambda_{1})=\mathbb{P}_{n}^{w}(\eta^{2})-\mathbb{P}_{n}^{w}(\eta)^{2},\\
&b_{1}^{n}(\lambda_{1})=\mathbb{P}_{n}^{w}(\eta Y^{2})-\mathbb{P}_{n}^{w}(\eta)\mathbb{P}_{n}^{w}(Y^{2}),\\
&b_{2}^{n}(\lambda_{1})=\mathbb{P}_{n}^{w}(\eta^{2}Y^{2})-\mathbb{P}_{n}^{w}(\eta^{2})\mathbb{P}_{n}^{w}(Y^{2}).
\end{align*}

Finally define

\begin{align*}
&\lambda_{2}^{n}(\lambda_{1})=\frac{a_{1}^{n}(\lambda_{1})\cdot b_{1}^{n}(\lambda_{1})-a_{2}^{n}(\lambda_{1})\cdot
b_{2}^{n}(\lambda_{1}) }{a_{1}^{n}(\lambda_{1})\cdot a_{3}^{n}(\lambda_{1}) - a_{2}^{n}(\lambda_{1})^{2}},\\
&\lambda_{3}^{n}(\lambda_{1})=\frac{a_{3}^{n}(\lambda_{1})\cdot b_{2}^{n}(\lambda_{1})-a_{2}^{n}(\lambda_{1})\cdot
b_{1}^{n}(\lambda_{1}) }{a_{1}^{n}(\lambda_{1})\cdot a_{3}^{n}(\lambda_{1}) - a_{2}^{n}(\lambda_{1})^{2}}.
\end{align*}

Then our estimators are given as

\begin{align}
&\hat{\lambda}_{2}=\lambda_{2}^{n}(\hat{\lambda}_{1}),\\
&\hat{\lambda}_{3}=\lambda_{3}^{n}(\hat{\lambda}_{1}).
\end{align}

We now proceed to derive large sample properties of the estimators. In
what follows we adopt the notation \(\mathbf{P}f=Ef(X_{i},Y_{i})\) for
any function of \(f\) of our data samples \((X_{i},Y_{i})\)

To derive large sample results we need to assume that the weight
function \(w\) is twice continuously differentiable and that
\(\lambda_{1}\rightarrow w(\lambda_{1}^{T}X)\) and the derivatives
\(\lambda_{1}\rightarrow D_{\lambda_{1}}w(\lambda_{1}^{T}X)\) and
\(\lambda_{1}\rightarrow D_{\lambda_{1}}^{2}w(\lambda_{1}^{T}X)\) are
bounded by functions with finite means in some open neighbourhood of the
true value of \(\lambda_{1}\).

We also need the following lemmas

\begin{lemma}\label{asympapprox}
Let $f:\mathbb{R}\times\mathbb{R}\rightarrow \mathbb{R}$ be a two times continuously differentiable  
function in the first coordinate. Assume that $\lambda_{1}\rightarrow f(\lambda_{1}^{T}X,Y)\cdot w(\lambda_{1}^{T}X)$ and the derivatives $\lambda_{1}\rightarrow D_{\lambda_{1}}\{f(\lambda_{1}^{T}X,Y)\cdot w(\lambda_{1}^{T}X)\}$ and $\lambda_{1}\rightarrow D^2_{\lambda_{1}}\{f(\lambda_{1}^{T}X,Y)\cdot w(\lambda_{1}^{T}X)\}$ are bounded by functions with finite mean in an open neighbourhood of the true value of $\lambda_{1}$. Also put $f_{i}=f(\eta_{i},Y_{i})$ and $\hat{f}_{i}=f(\hat{\eta}_{i},Y_{i})$. Then with 
$$
g(\lambda_{1})=\frac{\mathbf{P}\{f\cdot w\}}{\mathbf{P}w}
$$
it holds that $g$ is differentiable with derivative
$$
D_{\lambda_{1}}g(\lambda_{1})=(\mathbf{P}w)^{-1}P\{D_{\lambda_{1}}(f\cdot w)\}-(\mathbf{P}w)^{-2}\mathbf{P}\{f\cdot w\}\mathbf{P}(D_{\lambda_{1}}w).
$$
Moreover it holds that 
$$
\sqrt{n}\{\mathbb{P}_{n}^{\hat{w}}(\hat{f})-\mathbb{P}_{n}^{w}(f)\}=\sqrt{n}\{\hat{\lambda}_{1}-\lambda_{1}\}^{T}D_{\lambda_{1}}g(\lambda_{1})+o_{P}(1).
$$
\end{lemma}

\begin{lemma}\label{iidasymp}
Let $f:\mathbb{R}\times\mathbb{R}\rightarrow \mathbb{R}$ and put $f_{i}=f(\eta_{i},Y_{i})$.  Assume that $$
P\{f\cdot w\}^{2}<\infty.
$$
Then 

$$
\sqrt{n}\big\{\mathbb{P}^{w}_{n}(f)-\frac{\mathbf{P}\{f\cdot w\}}{\mathbf{P}w}\big\}=\{\mathbf{P}w\}^{-1}\mathbb{G}_{n}(\{f-\frac{\mathbf{P}(f\cdot w)}{\mathbf{P}w}\}\cdot w)+o_{P}(1),
$$
where $\mathbb{G}_{n}(f)=\frac{1}{\sqrt{n}}\sum_{i=1}^{n}(f_{i}-\mathbf{P}f)$.
\end{lemma}

\begin{theorem}\label{asympprop}

The estimators $\hat{\lambda}_{2}$ and $\hat{\lambda}_{3}$ are consistent, asymptotically normal, and have the following asymptotic characterizations: 

\begin{equation*}
\sqrt{n}(\hat{\lambda}_{2}-\lambda_{2})=\sqrt{n}\{\hat{\lambda}_{1}-\lambda_{1}\}^{T}D_{\lambda_{1}}\lambda_{2}(\lambda_{1})+\{g_{3}(\lambda_{1})\mathbf{P}w\}^{-1}\mathbf{G_{n}}(w\{a_{1}(\lambda_{1})f^{1}-a_{2}(\lambda_{1})f^{2}\})+o_{P}(1),
\end{equation*}

\begin{equation*}
\sqrt{n}(\hat{\lambda}_{3}-\lambda_{3})=\sqrt{n}\{\hat{\lambda}_{1}-\lambda_{1}\}^{T}D_{\lambda_{1}}\lambda_{3}(\lambda_{1})+\{g_{3}(\lambda_{1})\mathbf{P}w\}^{-1}\mathbf{G_{n}}(w\{a_{3}(\lambda_{1})f^{2}-a_{2}(\lambda_{1})f^{1}\})+o_{P}(1),
\end{equation*}
where $a_{1}(\lambda_{1})$, $a_{2}(\lambda_{1})$, $a_{3}(\lambda_{1})$, $g_{3}(\lambda_{1})$, $\lambda_{2}(\lambda_{1})$, and  $\lambda_{3}(\lambda_{1})$ are the limits in probability of $a^{n}_{1}(\lambda_{1})$, $a^{n}_{2}(\lambda_{1})$, $a^{n}_{3}(\lambda_{1})$, $g^{n}_{3}(\lambda_{1})$ $\lambda_{2}^{n}(\lambda_{1})$, and $\lambda_{3}^{n}(\lambda_{1})$, respectively, and

\begin{align*}
&f^{1}_{i}=(\eta_{i}-\frac{\mathbf{P}(\eta\cdot w)}{\mathbf{P}w})\{Y_{i}^{2}-E(Y_{i}^{2}|X_{i})\},\\
&f^{2}_{i}=(\eta_{i}^{2}-\frac{\mathbf{P}(\eta^{2}\cdot w)}{\mathbf{P}w})\{Y_{i}^{2}-E(Y_{i}^{2}|X_{i})\}.
\end{align*}
 
\end{theorem}

As an immediate consequence of (\ref{iid1}) and Theorem 1 we have the
that the obvious estimator
\(\hat{\beta}=\hat{\lambda}_{3}\hat{\lambda}_{1}\) of \(\beta\) is both
consistent and asymptotically normal. In particular, we have the
following iid decomposition of \(\sqrt{n}(\hat{\beta}-\beta)\):

\begin{equation}
\sqrt{n}(\hat{\beta}-\beta)=\frac{1}{\sqrt{n}}\sum_{i}^{n}\xi_{i}+o_{P}(1),
\end{equation}

where

\begin{equation*}
\xi_{i}=[\lambda_{3}\cdot I_{m}+\lambda_{1}\{D_{\lambda_{1}}\lambda_{3}(\lambda_{1})\}^{T}]\varepsilon_{i}+\lambda_{1}\cdot\{g_{3}(\lambda_{1}\mathbf{P}w\}^{-1}\cdot w_{i}\cdot\{a_{1}(\lambda_{1})f^{1}_{i}-a_{2}(\lambda_{1})f^{2}_{i}\}.
\end{equation*}

Notice that the iid decomposition above enables consistent estimation of
standard errors based on the law of large numbers. Specifically

\begin{equation}\label{seest}
\frac{1}{n}\sum_{i=1}^{n}\xi_{i}\xi_{i}^{T}\overset{P}\rightarrow \text{var}(\xi_{i})
\end{equation}

and accordingly a consistent estimator of \(\text{var}(\xi_{i})\) can be
obtained by inserting the empirical counterpart of \(\xi_{i}\) in the
left hand side of (\ref{seest}).

\subsection{Choice of weights}\label{choice-of-weights}

The optimal weights for weighted least squares regression are
\(w_i \propto \frac{1}{V(Y|X)}\). For the present we suggest using

\[v_i = v(\eta_i) = \frac{1}{(1 + \eta_i^2)}
\]

as weights for the first order moment regression, \eqref{firstorder},
and

\[w_i = w(\eta_i) = \frac{1}{(1 + \eta_i^4)}
\]

as weights for the second order moment regression \eqref{secorder}. Both
choices mimic the order of conditional variances as functions of
\(\eta_i\) and fulfills the requirements that the weight functions are
twice continuously differentiable and locally bounded by functions with
finite means. In the rest of the following we will use these weights.

\section{Simulation study}\label{simulation-study}

We illustrate the empirical properties of the proposed estimation
procedure through a series of simulations. Consider the following four
scenarios:

\begin{enumerate}
\def\labelenumi{\arabic{enumi}.}
\item
  Simple Gaussian mixture: \(Y_1 \sim N(1+X, 1^2), Y_2 \sim N(0,1^2)\)
  and \(X\sim N(0,1^2)\).
\item
  Gaussian distribution with zero-inflation:
  \(Y_1 \sim N(1+X, 1^2), Y_2 = 0\) and \(X\sim N(0,1^2)\).
\item
  Exponential-Gaussian mixture: \(Y_1 \sim \text{exp}(1)+X\) (shifted
  exponential distribution with rate \(1\)), \(Y_2 \sim N(0,0.5^2)\) and
  \(X\sim N(0,1^2)\).
\item
  Exponential distribution with zero-inflation:
  \(Y_1 \sim \text{exp}(1)+X\), \(Y_2=0\) and \(X\sim N(0,1^2)\).
\end{enumerate}

For each scenario we simulate 100 datasets each with
\(N\in\{300, 500, 800, 100, 1500, 2000, 3000 \}\) observations and with
a constant mixture proportion of \(p=0.7\) (i.e., 30\% contamination).
The estimates \(\hat\beta\) and \(\hat p\) are computed using both the
moment mixture estimation procedure as well as a standard two-component
Gaussian regression mixture model where \(X\) can only influence one of
the components.

\begin{figure}
\centering
\includegraphics{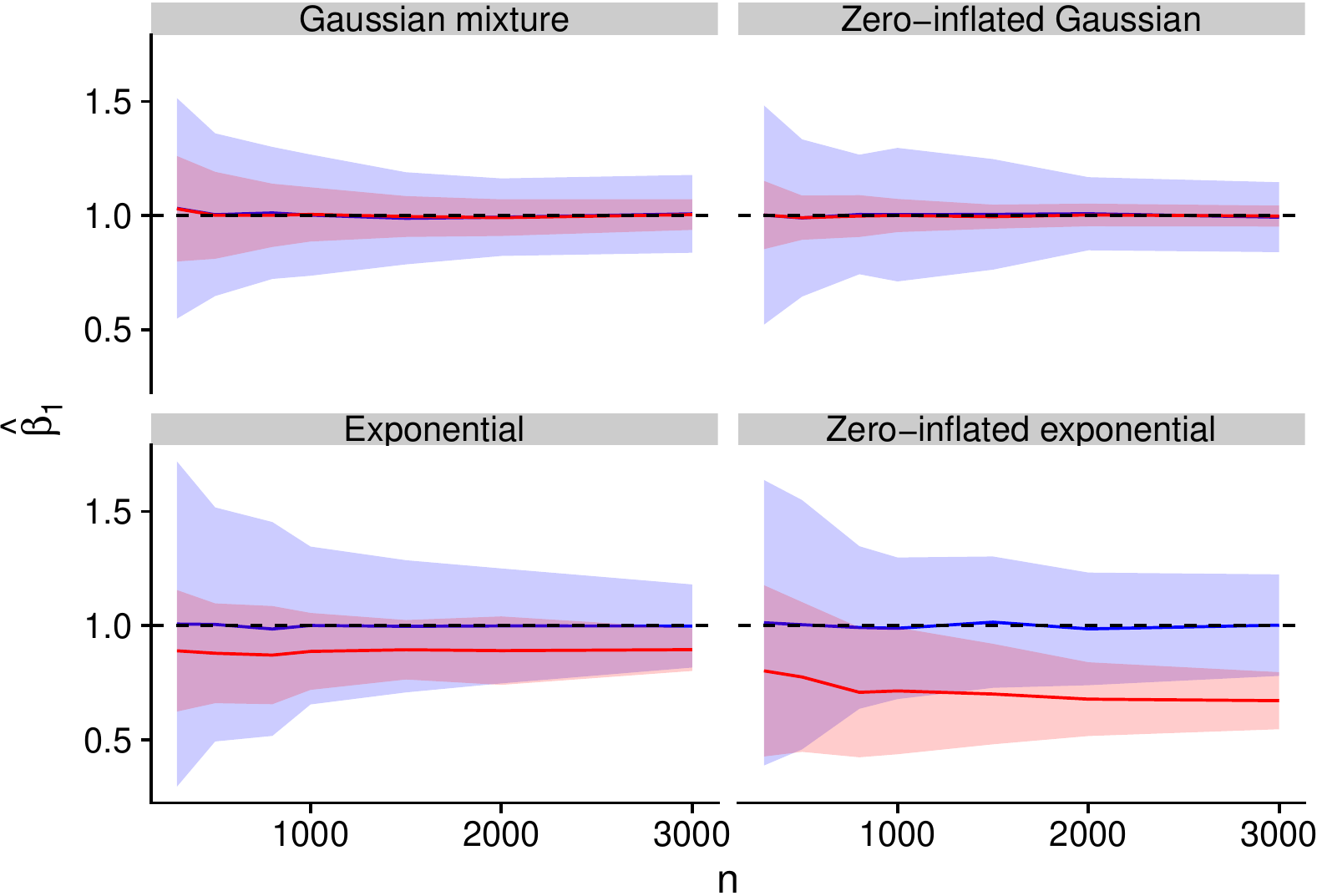}
\caption{\label{fig:fig1}Average estimates of \(\hat\beta\) using the moment
mixture model (blue) and Gaussian mixture model (red) and corresponding
95\% confidence pointwise confidence bands.}
\end{figure}

\begin{figure}
\centering
\includegraphics{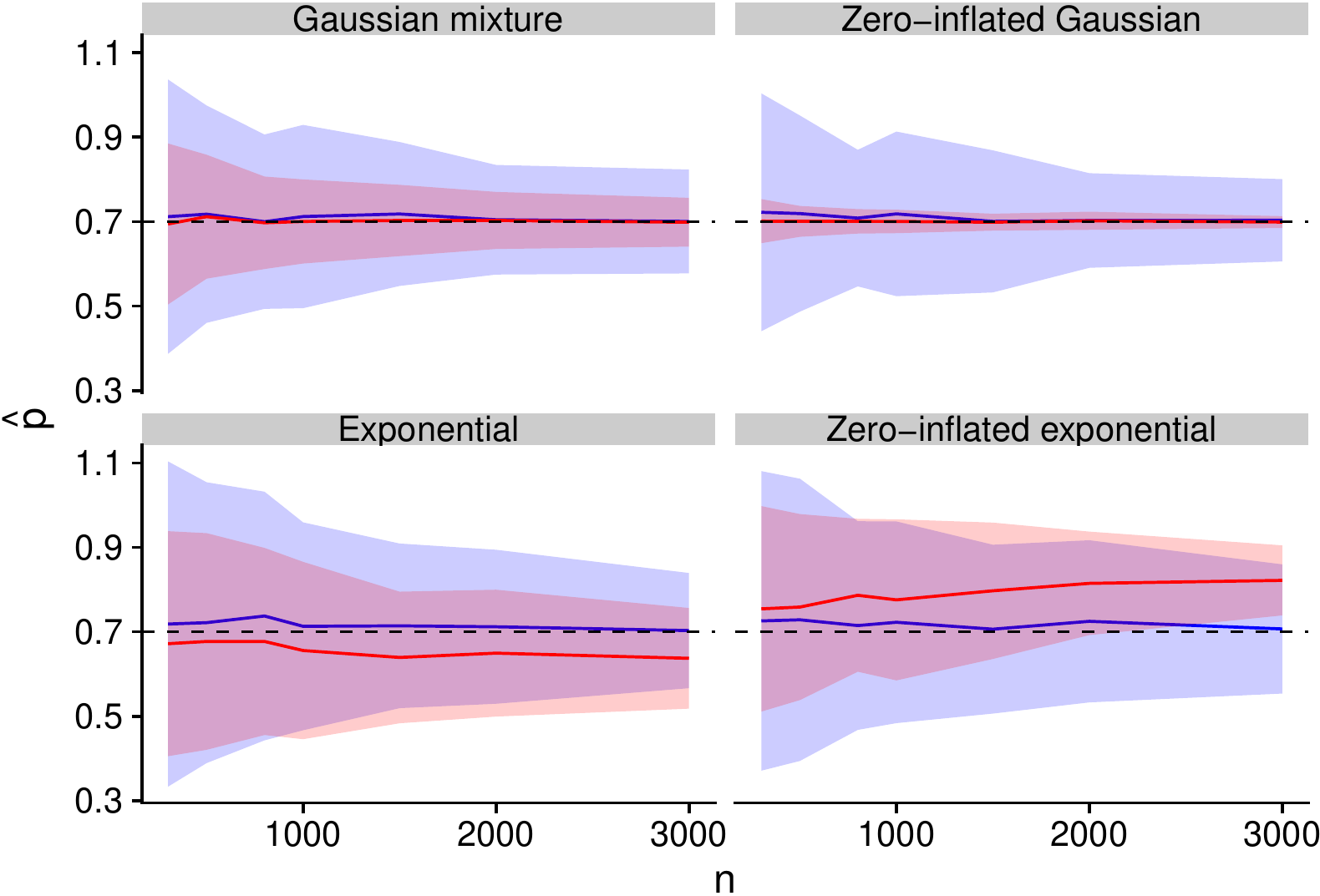}
\caption{\label{fig:fig2}Average estimates of mixing proportion \(\hat p\)
using the moment mixture model (blue) and Gaussian mixture model (red)
and corresponding 95\% confidence pointwise confidence bands.}
\end{figure}

Figure \ref{fig:fig1} shows the empirical estimates of \(\hat\beta\)
(and the corresponding 95\% pointwise standard errors) for the four
different scenarios while Figure \ref{fig:fig1} shows the corresponding
plot for the mixing proportion \(\hat p\). Not surprisingly, the
Gaussian regression mixture model provides an unbiased estimate of both
\(\hat\beta\) and \(\hat p\) with the smallest variance when the data
follows a mixture of Gaussian distributions. However, it is also clear
that the moment mixture model --- which places fewer assumptions on the
distributions involved --- is also unbiased and has a variance that is
only slightly larger than the Gaussian mixture model. Substantial
differences are seen when the distribution of \(Y_1\) is \emph{not}
Gaussian where the moment mixture model quickly converges to the true
value while the Gaussian mixture model results in biased estimates of
both \(\beta\) and \(p\). The estimates from the moment mixture model
are unbiased even when the distributions are non-normal, but Figure
\ref{fig:fig2} suggests that large samples may be needed in order to
have sufficient power to use the moment mixture model to infer whether
there are two or only one component.

The variance of the mixing proportion estimate, \(\text{var}(\hat p)\),
is larger for the moment mixture model than for the Gaussian regression
mixture model which essentially is the ``price that is paid'' by having
fewer assumptions about the distribution and estimating the parameters
from the moments of the distribution.

Table \ref{tab:tab1} shows results from the same four scenarios
(although with mixing proportion \(p=0.5\)) to illustrate the precision
of the estimates and to compare the estimated standard errors from
\eqref{seest} to the empirical standard errors. The estimates for both
\(\beta\) and \(p\) are close to the true values even for smaller sample
sizes and when the error distribution is non-normal. Table
\ref{tab:tab1} also shows that the estimated standard errors are close
to the empirical standard errors which suggests that the asymptotic
estimate in \eqref{seest} provides a useful measure of the standard
errors of the regression and mixing proportion parameters even for
smaller sample sizes. The 95\% coverage probabilities for \(\beta\) are
close to the true value even for \(N=300\) regardless of the underlying
distributions while the coverage probabilities for \(p\) are quite
unstable for smaller sample sizes while they achieve right right level
for \(N=2000\). Note that the mixing proportion is 50\% so there is a
substantial amount of noise in the data.

To investigate the efficiency of the proposed estimates we consider the
case where we know the individual allocations, \(P_{i}\). If we make no
assumptions on the distribution of \(Y_{2}\) but assume that \(Y_{1}\)
follows a normal distribution with mean \(\beta_{0}+\beta_{1} X\) and
variance \(\sigma_{1}^{2}\) conditional on \(X\in\mathbb{R}\), then the
most efficient way to estimate \(\beta_{1}\) would be by simple linear
regression based only on the pairs \((X_{i},Y_{i})\) where \(P_{i}=1\).
For this estimator the asymptotic variance is equal to

\begin{equation}\label{effbound}
\text{Var}(\hat\beta) = \frac{\sigma_{1}^{2}}{(p\cdot N -1) \cdot \text{Var}(X_{i})}
\end{equation}

Accordingly, the maximum likelihood estimator \(\hat{\beta}_{1}\) of
\(\beta_{1}\) from any parametric mixture model of \(Y\) with the above
specification of \(Y_{1}\) will have an asymptotic variance that is at
least (\ref{effbound}). We compare the estimated standard error from the
moment mixture model to this lower bound to see the efficiency loss
incurred by removing the assumptions about the parametric distributions.
The result is seen in the right-most column of Table \ref{tab:tab1}, and
we see that the efficiency loss from using the moment mixture model is
around a factor 4-5.

\rowcolors{2}{gray!6}{white}

\begin{table}

\caption{\label{tab:tab1}Results from simulation study of the moment mixture model based on 1000 simulations. The model provides unbiased estimates and the estimted SE match the empirical SEs. 95\% CP is the 95\% coverage probability and the efficiency is the relative efficiency of the SE for the moment mixture model relative to the SE of the optimal model.}
\centering
\begin{tabular}[t]{lrrrrrrr}
\hiderowcolors
\toprule
  & N & True value & Estimate & Emp. SE & Est. SE & 95\% CP & Rel. efficiency\\
\midrule
\showrowcolors
\addlinespace[0.3em]
\multicolumn{8}{l}{\textbf{Gaussian mixture}}\\
\hspace{1em}$\beta$ & 300 & 1.0 & 1.004 & 0.331 & 0.323 & 0.939 & 3.96\\
\hspace{1em}p & 300 & 0.5 & 0.481 & 0.197 & 0.174 & 0.850 & \\
\hspace{1em}$\beta$ & 2000 & 1.0 & 1.002 & 0.120 & 0.121 & 0.948 & 3.83\\
\hspace{1em}p & 2000 & 0.5 & 0.497 & 0.066 & 0.065 & 0.948 & \\
\addlinespace[0.3em]
\multicolumn{8}{l}{\textbf{Zero-inflated Gaussian}}\\
\hspace{1em}$\beta$ & 300 & 1.0 & 0.977 & 0.315 & 0.307 & 0.943 & 3.76\\
\hspace{1em}p & 300 & 0.5 & 0.506 & 0.181 & 0.156 & 0.870 & \\
\hspace{1em}$\beta$ & 2000 & 1.0 & 0.995 & 0.117 & 0.118 & 0.956 & 3.73\\
\hspace{1em}p & 2000 & 0.5 & 0.502 & 0.059 & 0.056 & 0.957 & \\
\addlinespace[0.3em]
\multicolumn{8}{l}{\textbf{Exponential-Gaussian mixture}}\\
\hspace{1em}$\beta$ & 300 & 1.0 & 0.998 & 0.441 & 0.394 & 0.929 & 4.83\\
\hspace{1em}p & 300 & 0.5 & 0.476 & 0.231 & 0.204 & 0.712 & \\
\hspace{1em}$\beta$ & 2000 & 1.0 & 1.006 & 0.153 & 0.154 & 0.953 & 4.87\\
\hspace{1em}p & 2000 & 0.5 & 0.495 & 0.085 & 0.080 & 0.956 & \\
\addlinespace[0.3em]
\multicolumn{8}{l}{\textbf{Zero-inflated exponential}}\\
\hspace{1em}$\beta$ & 300 & 1.0 & 0.996 & 0.404 & 0.375 & 0.942 & 4.59\\
\hspace{1em}p & 300 & 0.5 & 0.490 & 0.212 & 0.185 & 0.773 & \\
\hspace{1em}$\beta$ & 2000 & 1.0 & 0.999 & 0.154 & 0.151 & 0.946 & 4.78\\
\hspace{1em}p & 2000 & 0.5 & 0.500 & 0.080 & 0.073 & 0.947 & \\
\bottomrule
\end{tabular}
\end{table}

\rowcolors{2}{white}{white}

\section{Application: The effect of pH on wine volatile
acidity}\label{application-the-effect-of-ph-on-wine-volatile-acidity}

Volatile acidity (VA) refers to the process when lactic acid bacteria
and acetic acids turns wine into vinegar, and the process takes place
mainly due to growth of bacteria, the oxidation of ethanol, or the
metabolism of acids/sugars. Wines with a high level of pH are supposedly
more susceptible to oxidation and the antibacterial effects of sulfur
dioxide and of fumaric acid are reduced rapidly as the pH level
increases. Consequently, wines are thought to lose their quality as they
become less acidic (increased pH) since the volatile acidity increases.

\begin{figure}
\centering
\includegraphics{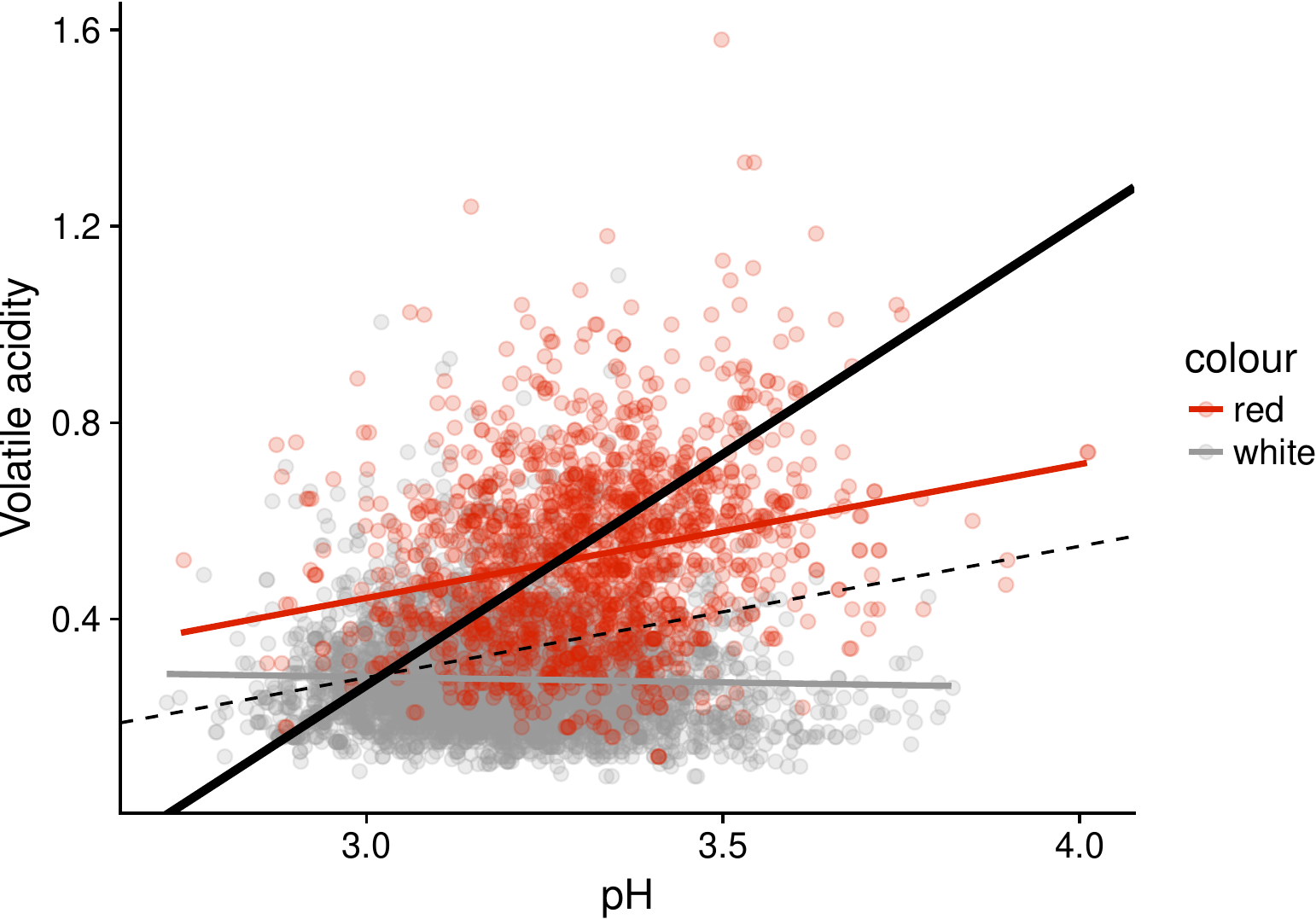}
\caption{\label{fig:wine}Volatile acidity vs.~pH for white and red wines.
The coloured lines show the estimated regression lines for the two
corresponding types of wines, respectively. The dashed line shows the
estimated regression line when a linear regression is used on all data
points, while the solid line shows the estimated regression line for the
moment mixture model.}
\end{figure}

The paper by \citet{cort:etal:2009} considers 11 physicochemical
properties of a selected sample of Portuguese \emph{vinho verde} wines.
Samples from 1599 red wines and 4898 white wines are available and the
relationship between volatile acidity and pH are shown separately for
red and white wines in Figure \ref{fig:wine}. While it is apparent that
the available red wines generally have slightly higher levels of pH, it
also appears as if the impact of pH on volatile acidity is largest for
the red wines in the sample: Individual regression lines for the two
types of wine show an almost horizontal line for white wine while there
is an effect of pH on VA for the red wines (slopes \(-0.022\) and
\(0.272\), respectively, and the slopes are significantly different,
\(p<0.0001\)).

If we did not know that the full 6497 samples were comprised of two
different types of wine we might pursue regressing volatile acidity on
pH for the full data. This gives the dashed regression shown in Figure
\ref{fig:wine} (slope \(\hat\beta\) =0.27, 95\% CI 0.24-0.29), which
suggests an overall effect of pH on volatile acidity.

When we fit the moment-based mixture model then we  see that the
data are likely to consist of a mixture, \(\hat{p} = 0.28\) (95\% CI:
0.24-0.33), which suggests that only a proportion of the wines are
influenced by pH. Also, the effect driving the association between the
volatile acidity and pH appears to be much stronger than what was
observed from analyzing the full data with a simple linear regression
model, \(\hat\beta=\) 0.94 (95\% CI 0.51-1.37). Consequently, the moment
mixture model is able to identify that the wine data is likely to
consist of two types of wines that respond differently to changes in
levels of pH with only placing very few restrictions on the underlying
distributions.

We estimate the mixing proportion to be 28\% whereas the dataset
contains 24.6\% red wines. While the wine colour classification need not
correspond to the separation we estimate from the moment mixture model
it may indeed be the case that the different wine types contain a set of
features that influence the impact of pH which is what we observe for
these data.

\section{Discussion}\label{discussion}

In this article we address the problem of estimating regression
parameters for a two-component mixture regression model where one
component is influenced by a set of predictors. Our proposed method
imposes no assumptions on the distributions of the components and only a
minimum of restrictions on the regression effects. The moment-based
mixture of regression models estimators can be used to detect or account
for unspecified mixtures in regression problems, and since estimation is
fast it could be used for large-scale studies such as in genome-wise
association studies.

As shown by the simulation study, the price for the flexibility and lack
of assumptions comes with a larger variance of the estimators but the
variance is not prohibitively larger and the proposed methods provides
consistent estimates even when the mixture distributions are highly
skewed and other models such as the Gaussian mixture model fails.

In conclusion, we have introduced an estimation technique for mixtures
of regression models that can be applied to a large number of
situations. The moment-based estimator is very versatile: it can be used
not only to estimate the regression parameters, but for larger datasets
it provides a foundation for detecting the number of mixture components.
If \(\hat p\) is different from both 0 and 1 then this suggests that
there indeed are two components where only one of them is influenced by
the predictors of interest.

\section{Software}
\label{sec5}

Software in the form of the \texttt{mommix} R package and code, together with a sample
input data set and complete documentation is available on
GitHub at \texttt{www.github.com/ekstroem/mommix}.


{\it Conflict of Interest}: None declared.

\bibliographystyle{biorefs}
\bibliography{refs}

\begin{thebibliography}{99}

\bibitem[Cortez \emph{and others}(2009)Cortez, Cerdeira, Almeida, Matos and
  Reis]{cort:etal:2009}
\textsc{Cortez, Paulo, Cerdeira, Antonio, Almeida, Fernando, Matos, Telmo and
  Reis, Jose}. (2009).
\newblock Modeling wine preferences by data mining from physicochemical
  properties.
\newblock {\em Decision Support Systems\/}~\textbf{47}, 557--53.

\bibitem[Dalrymple \emph{and others}(2003)Dalrymple, Hudson and
  Ford]{dalr:etal:2003}
\textsc{Dalrymple, M.L., Hudson, I.L. and Ford, R.P.K.} (2003).
\newblock Finite mixture, zero-inflated poisson and hurdle models with
  application to sids.
\newblock {\em Computational Statistics \& Data Analysis\/}~\textbf{41}(3), 491
  -- 504.
\newblock Recent Developments in Mixture Model.

\bibitem[Deb and Trivedi(1997)Deb and Trivedi]{deb:triv:1997}
\textsc{Deb, Partha and Trivedi, Pravin~K.} (1997).
\newblock Demand for medical care by the elderly: a finite mixture approach.
\newblock {\em Journal of Applied Econometrics\/}~\textbf{12}, 313--336.

\bibitem[Hu \emph{and others}(2011)Hu, Pavlicova and Nunes]{hu:etal:2011}
\textsc{Hu, Mei-Chen, Pavlicova, Martina and Nunes, Edward~V.} (2011).
\newblock Zero-inflated and hurdle models of count data with extra zeros:
  Examples from an {HIV}-risk reduction intervention trial.
\newblock {\em Am J Drug Alcohol Abuse\/}~\textbf{37}, 367--375.

\bibitem[Kong \emph{and others}(2017)Kong, Li, Fan and Lv]{kong:etal:2017}
\textsc{Kong, Yinfei, Li, Daoji, Fan, Yingying and Lv, Jinchi}. (2017).
\newblock Interaction pursuit in high-dimensional multi-response regression via
  distance correlation.
\newblock {\em Annals of Statistics\/}~\textbf{45}(2), 897--922.

\bibitem[Kwagyan and Apprey(2016)Kwagyan and Apprey]{kwag:appr:2016}
\textsc{Kwagyan, John and Apprey, Victor}. (2016).
\newblock Modeling clustered binary data with excess zero clusters.
\newblock {\em Statistical Methods in Medical Research\/}~\textbf{0}, 1--16.

\bibitem[Lambert(1992)Lambert]{lamb:1992}
\textsc{Lambert, D.} (1992).
\newblock Zero-inflated poisson regression with an application to defects in
  manufacturing.
\newblock {\em Technometrics\/}~\textbf{34}, 1--14.

\bibitem[McLachlan and Peel(2000)McLachlan and Peel]{mcla:peel:2000}
\textsc{McLachlan, G.~J. and Peel, D.} (2000).
\newblock {\em Finite mixture models\/}. New York: Wiley Series in Probability
  and Statistics.

\bibitem[Xu \emph{and others}(2015)Xu, Paterson, Turpin and Xu]{xu:etal:2015}
\textsc{Xu, Lizhen, Paterson, Andrew~D., Turpin, Williams and Xu, Wei}. (2015,
  07).
\newblock Assessment and selection of competing models for zero-inflated
  microbiome data.
\newblock {\em PLOS ONE\/}~\textbf{10}(7), 1--30.

\end{thebibliography}

\section{Appendix}\label{appendix}

\subsection{Proof of Theorem \ref{asympprop}}\label{proof-of-theorem}

Denote

\begin{align*}
&g_{1}^{n}(\lambda_{1})=a_{1}^{n}(\lambda_{1})\cdot b_{1}^{n}(\lambda_{1})-a_{2}^{n}(\lambda_{1})\cdot
b_{2}^{n}(\lambda_{1}),\\
&g^{n}_{2}(\lambda_{1})=a_{3}^{n}(\lambda_{1})\cdot b_{2}^{n}(\lambda_{1})-a_{2}^{n}(\lambda_{1})\cdot
b_{1}^{n}(\lambda_{1}),\\
&g^{n}_{3}(\lambda_{1})=a_{1}^{n}(\lambda_{1})\cdot a_{3}^{n}(\lambda_{1}) - a_{2}^{n}(\lambda_{1})^{2}
\end{align*}

with corresponding limits in probability

\begin{align*}
&g_{1}(\lambda_{1})=a_{1}(\lambda_{1})\cdot b_{1}(\lambda_{1})-a_{2}(\lambda_{1})\cdot
b_{2}(\lambda_{1}),\\
&g_{2}(\lambda_{1})=a_{3}(\lambda_{1})\cdot b_{2}(\lambda_{1})-a_{2}(\lambda_{1})\cdot
b_{1}(\lambda_{1}),\\
&g_{3}(\lambda_{1})=a_{1}(\lambda_{1})\cdot a_{3}(\lambda_{1}) - a_{2}(\lambda_{1})^{2},
\end{align*}

where

\begin{align*}
&a_{1}(\lambda_{1})=\frac{\mathbf{P}(\eta^{4}\cdot w)}{\mathbf{P} w}-\big[\frac{\mathbf{P}(\eta^{2}\cdot w)}{\mathbf{P} w}\big]^{2},\\
&a_{2}(\lambda_{1})=\frac{\mathbf{P}(\eta^{3}\cdot w)}{\mathbf{P} w}-\frac{\mathbf{P}(\eta \cdot w)\cdot\mathbf{P}(\eta^{2}\cdot w)}{(\mathbf{P} w)^{2}},\\
&a_{3}(\lambda_{1})=\frac{\mathbf{P}(\eta^{2}\cdot w)}{\mathbf{P}w}-\big[\frac{\mathbf{P}(\eta\cdot w)}{\mathbf{P}w}\big]^{2},\\
&b_{1}(\lambda_{1})=\frac{\mathbf{P}(\eta\cdot Y^{2}\cdot w)}{\mathbf{P}w}-\frac{\mathbf{P}(\eta\cdot w)\mathbf{P}(Y^{2}\cdot w)}{(\mathbf{P}w)^{2}},\\
&b_{2}(\lambda_{1})=\frac{\mathbf{P}(\eta^{2}\cdot Y^{2}\cdot w)}{\mathbf{P}w}-\frac{\mathbf{P}(\eta^{2}\cdot w)\mathbf{P}(Y^{2}\cdot w)}{(\mathbf{P}w)^{2}}.
\end{align*}

By repeatedly using Lemma \ref{asympapprox}, the consistency of
\(\hat{\lambda}_{1}\), and standard arguments we have

\begin{align*}
&\sqrt{n}\{g_{1}^{n}(\hat{\lambda}_{1})-g^{1}_{n}(\lambda_{1})\}=\sqrt{n}\{\hat{\lambda}_{1}-\lambda_{1}\}^{T}D_{\lambda_{1}}g_{1}(\lambda_{1})+o_{P}(1),\\
&\sqrt{n}\{g_{2}^{n}(\hat{\lambda}_{1})-g_{2}^{n}(\lambda_{1})\}=\sqrt{n}\{\hat{\lambda}_{1}-\lambda_{1}\}^{T}D_{\lambda_{1}}g_{2}(\lambda_{1})+o_{P}(1),\\
&\sqrt{n}\{g_{3}^{n}(\hat{\lambda}_{1})-g_{3}^{n}(\lambda_{1})\}=\sqrt{n}\{\hat{\lambda}_{1}-\lambda_{1}\}^{T}D_{\lambda_{1}}g_{3}(\lambda_{1})+o_{P}(1).
\end{align*}

Then\(\sqrt{n}(\hat{\lambda}_{2}-\lambda_{2})\) may be rewritten as:

\begin{equation}\label{first}
\sqrt{n}(\hat{\lambda}_{2}-\lambda_{2})=\sqrt{n}\Big(\frac{g^{n}_{1}(\hat{\lambda}_{1})}{g^{n}_{3}(\hat{\lambda}_{1})}-\frac{g^{n}_{1}(\lambda_{1})}{g^{n}_{3}(\lambda_{1})}\Big)+
\sqrt{n}\Big(\frac{g^{n}_{1}(\lambda_{1})}{g^{n}_{3}(\lambda_{1})}-\lambda_{2}\Big)
\end{equation}

For the first term of the right hand side of equation (\ref{first})

\begin{multline}
  \sqrt{n}\Big(\frac{g^{n}_{1}(\hat{\lambda}_{1})}{g^{n}_{3}(\hat{\lambda}_{1})}-\frac{g^{n}_{1}(\lambda_{1})}{g^{n}_{3}(\lambda_{1})}\Big) \\
  =  [g_{3}(\lambda_{1})]^{-1}\sqrt{n}\{g^{n}_{1}(\hat{\lambda}_{1})-g^{n}_{1}(\lambda_{1})\}-[g_{3}(\lambda_{1})]^{-2}[g_{1}(\lambda_{1})]\sqrt{n}\{g_{3}^{n}(\hat{\lambda}_{1})-g_{3}^{n}(\lambda_{1}))\}+o_{P}(1)
\end{multline}

by the consistency of \(\hat{\lambda}_{1}\) and the law of large
numbers. Combining the above results we obtain

\begin{equation}\label{second}
\sqrt{n}\Big(\frac{g^{n}_{1}(\hat{\lambda}_{1})}{g^{n}_{3}(\hat{\lambda}_{1})}-\frac{g^{n}_{1}(\lambda_{1})}{g^{n}_{3}(\lambda_{1})}\Big)=\sqrt{n}\{\hat{\lambda}_{1}-\lambda_{1}\}^{T}D_{\lambda_{1}}\{\frac{g_{1}(\lambda_{1})}{g_{3}(\lambda_{1})}\}+o_{P}(1)=\sqrt{n}\{\hat{\lambda}_{1}-\lambda_{1}\}^{T}D_{\lambda_{1}}\lambda_{2}(\lambda_{1})+o_{P}(1).
\end{equation}

For the second term on the right hand side of equation (\ref{first})
first note that by repeatedly using Lemma \ref{iidasymp} and the law of
large numbers

\begin{multline}
\sqrt{n}\{a^{n}_{1}(\lambda_{1})-a_{1}(\lambda_{1})\}=\\
(\mathbf{P}w)^{-1}\cdot\mathbb{G}_{n}(\{\eta^{4}-2\frac{\mathbf{P}(\eta^{2}w)}{\mathbf{P}w}\eta^{2}-a_{1}(\lambda_{1})+[\frac{\mathbf{P}(\eta^{2}w)}{\mathbf{P}w}]^{2}\}w)+o_{P}(1),\\
\sqrt{n}\{a^{n}_{2}(\lambda_{1})-a_{2}(\lambda_{1})\}= \\
(\mathbf{P}w)^{-1}\cdot\mathbb{G}_{n}(\{\eta^{3}-\frac{\mathbf{P}(\eta\cdot w)}{\mathbf{P}w}\eta^{2}-\frac{\mathbf{P}(\eta^{2}w)}{\mathbf{P}w}\eta-a_{2}(\lambda_{1})+\frac{\mathbf{P}(\eta \cdot w)\cdot\mathbf{P}(\eta^{2}\cdot w)}{(\mathbf{P} w)^{2}}\}w)+o_{P}(1),\\
\sqrt{n}\{a^{n}_{3}(\lambda_{1})-a_{3}(\lambda_{1})\}=\\(\mathbf{P}w)^{-1}\cdot\mathbb{G}_{n}(\{\eta^{2}-2\frac{\mathbf{P}(\eta\cdot w)}{\mathbf{P}w}\eta-a_{3}(\lambda_{1})+[\frac{\mathbf{P}(\eta\cdot w)}{\mathbf{P}w}]^{2}\}w)
+o_{P}(1),\\
\sqrt{n}\{b^{n}_{1}(\lambda_{1})-b_{1}(\lambda_{1})\}=\\
(\mathbf{P}w)^{-1}\cdot\mathbb{G}_{n}(\{\eta\cdot Y^{2}-\frac{\mathbf{P}(\eta\cdot w)}{\mathbf{P}w}Y^{2}-\frac{\mathbf{P}(Y^{2}w)}{\mathbf{P}w}\eta-b_{1}(\lambda_{1})+\frac{\mathbf{P}(\eta\cdot w)\mathbf{P}(Y^{2}\cdot w)}{(\mathbf{P}w)^{2}}\}w)
+o_{P}(1),\\
\sqrt{n}\{b^{n}_{2}(\lambda_{1})-b_{2}(\lambda_{1})\}=\\
(\mathbf{P}w)^{-1}\cdot\mathbb{G}_{n}(\{\eta^{2} Y^{2}-\frac{\mathbf{P}(\eta^{2} w)}{\mathbf{P}w}Y^{2}-\frac{\mathbf{P}(Y^{2}w)}{\mathbf{P}w}\eta^{2}-b_{2}(\lambda_{1})+\frac{\mathbf{P}(\eta^{2} w)\mathbf{P}(Y^{2}\cdot w)}{(\mathbf{P}w)^{2}}\}w)
+o_{P}(1).
\end{multline}

Next, by the law of large numbers

\begin{align*}
\sqrt{n}\{g_{1}^{n}(\lambda_{1})-g_{1}(\lambda_{1})\}=&a_{1}(\lambda_{1})\sqrt{n}\{b^{n}_{1}(\lambda_{1})-b_{1}(\lambda_{1})\}-a_{2}(\lambda_{1})\sqrt{n}\{b^{n}_{2}(\lambda_{1})-b_{2}(\lambda_{1})\}\\
&+b_{1}(\lambda_{1})\sqrt{n}\{a^{n}_{1}(\lambda_{1})-a_{1}(\lambda_{1})\}-b_{2}(\lambda_{1})\sqrt{n}\{a^{n}_{2}(\lambda_{1})-a_{2}(\lambda_{1})\}+o_{P}(1),\\
\sqrt{n}\{g_{2}^{n}(\lambda_{1})-g_{2}(\lambda_{1})\}=&-a_{2}(\lambda_{1})\sqrt{n}\{b^{n}_{1}(\lambda_{1})-b_{1}(\lambda_{1})\}+a_{3}(\lambda_{1})\sqrt{n}\{b^{n}_{2}(\lambda_{1})-b_{2}(\lambda_{1})\}\\
&-b_{1}(\lambda_{1})\sqrt{n}\{a^{n}_{2}(\lambda_{1})-a_{2}(\lambda_{1})\}+b_{2}(\lambda_{1})\sqrt{n}\{a^{n}_{3}(\lambda_{1})-a_{3}(\lambda_{1})\}+o_{P}(1),\\
\sqrt{n}\{g_{3}^{n}(\lambda_{1})-g_{3}(\lambda_{1})\}=&a_{3}(\lambda_{1})\sqrt{n}\{a_{1}^{n}(\lambda_{1})-a_{1}(\lambda_{1})\}+a_{1}(\lambda_{1})\sqrt{n}\{a_{3}^{n}(\lambda_{1})-a_{3}(\lambda_{1})\}\\
&-2\cdot a_{2}(\lambda_{1})\sqrt{n}\{a_{2}^{n}(\lambda_{1})-a_{2}(\lambda_{1})\}+o_{P}(1).
\end{align*}

Now note that in the true value of \(\lambda\) we have

\begin{align*}
&b_{1}(\lambda_{1})=\lambda_{2}a_{3}(\lambda_{1})+\lambda_{3}a_{2}(\lambda_{1}),\\
&b_{2}(\lambda_{1})=\lambda_{2}a_{2}(\lambda_{1})+\lambda_{3}a_{1}(\lambda_{1})
\end{align*}

and thus

\begin{align*}
&g_{1}(\lambda_{1})-\lambda_{2}g_{3}(\lambda_{1})=0,\\
&g_{2}(\lambda_{1})-\lambda_{3}g_{3}(\lambda_{1})=0.
\end{align*}

Accordingly

\begin{align*}
&\sqrt{n}\{g_{1}^{n}(\lambda_{1})-\lambda_{2}g_{3}^{n}(\lambda_{1})\}=\sqrt{n}\{g_{1}^{n}(\lambda_{1})-g_{1}(\lambda_{1})\}-\lambda_{2}\sqrt{n}\{g_{3}^{n}(\lambda_{1})-g_{3}(\lambda_{1})\},\\
&\sqrt{n}\{g_{2}^{n}(\lambda_{1})-\lambda_{3}g_{3}^{n}(\lambda_{1})\}=\sqrt{n}\{g_{2}^{n}(\lambda_{1})-g_{2}(\lambda_{1})\}-\lambda_{3}\sqrt{n}\{g_{3}^{n}(\lambda_{1})-g_{3}(\lambda_{1})\}.
\end{align*}

Finally, combining all of the above and using the law of large numbers:

\begin{equation}\label{third}
\sqrt{n}\Big(\frac{g^{n}_{1}(\lambda_{1})}{g^{n}_{3}(\lambda_{1})}-\lambda_{2}\Big)=\{g_{3}(\lambda_{1})\mathbf{P}w\}^{-1}\mathbf{G_{n}}(w\{a_{1}(\lambda_{1})f^{1}-a_{2}(\lambda_{1})f^{2}\}),
\end{equation}

where

\begin{align*}
f^{1}_{i}=&\eta_{i}Y_{i}^{2}-\frac{\mathbf{P}(\eta\cdot w)}{\mathbf{P}w}Y_{i}^{2}-\lambda_{3}\eta_{i}^{3}-\{\lambda_{2}-\lambda_{3}\frac{\mathbf{P}(\eta\cdot w)}{\mathbf{P}w}\}\eta_{i}^{2}-\\
& \{\frac{\mathbf{P}(Y^{2}w)}{\mathbf{P}w}-\lambda_{2}\frac{2\mathbf{P}(\eta\cdot w)}{\mathbf{P}w}-\lambda_{3}\frac{\mathbf{P}(\eta^{2}w)}{\mathbf{P}w}\}\eta_{i}\\
&+\frac{\mathbf{P}(\eta\cdot w)}{\mathbf{P}w}\{\frac{\mathbf{P}(Y^{2}w)}{\mathbf{P}w}-\lambda_{2}\frac{\mathbf{P}(\eta\cdot w)}{\mathbf{P}w}-\lambda_{3}\frac{\mathbf{P}(\eta^{2}w)}{\mathbf{P}w}\}\\
=&(\eta_{i}-\frac{\mathbf{P}(\eta\cdot w)}{\mathbf{P}w})\{Y_{i}^{2}-E(Y_{i}^{2}|X_{i})\},\\
f^{2}_{i}=&\eta_{i}^{2}Y_{i}^{2}-\frac{\mathbf{P}(\eta^{2}\cdot w)}{\mathbf{P}w}Y_{i}^{2}-\lambda_{3}\eta_{i}^{4}-\lambda_{2}\eta_{i}^{3}-\{\frac{\mathbf{P}(Y^{2}w)}{\mathbf{P}w}-\lambda_{3}\frac{2\mathbf{P}(\eta^{2}\cdot w)}{\mathbf{P}w}-\lambda_{2}\frac{\mathbf{P}(\eta \cdot w)}{\mathbf{P}w}\}\eta_{i}^{2}\\
&+\lambda_{2}\frac{\mathbf{P}(\eta^{2}w)}{\mathbf{P}w}\eta_{i}+\frac{\mathbf{P}(\eta^{2}\cdot w)}{\mathbf{P}w}\{\frac{\mathbf{P}(Y^{2}w)}{\mathbf{P}w}-\lambda_{3}\frac{\mathbf{P}(\eta^{2}\cdot w)}{\mathbf{P}w}-\lambda_{2}\frac{\mathbf{P}(\eta\cdot w)}{\mathbf{P}w}\}\\
=&(\eta_{i}^{2}-\frac{\mathbf{P}(\eta^{2}\cdot w)}{\mathbf{P}w})\{Y_{i}^{2}-E(Y_{i}^{2}|X_{i})\}.
\end{align*}

Adding (\ref{second}) and (\ref{third}) we obtain

\[
\sqrt{n}\{\hat{\lambda}_{2}-\lambda_{2}\}=\sqrt{n}\{\hat{\lambda}_{1}-\lambda_{1}\}^{T}D_{\lambda_{1}}\lambda_{2}(\lambda_{1})+\{g_{3}(\lambda_{1})\mathbf{P}w\}^{-1}\mathbf{G_{n}}(w\{a_{1}(\lambda_{1})f^{1}-a_{2}(\lambda_{1})f^{2}\})+o_{P}(1)
\] Similarly one may show \[
\sqrt{n}\{\hat{\lambda}_{3}-\lambda_{3}\}=\sqrt{n}\{\hat{\lambda}_{1}-\lambda_{1}\}^{T}D_{\lambda_{1}}\lambda_{3}(\lambda_{1})+\{g_{3}(\lambda_{1})\mathbf{P}w\}^{-1}\mathbf{G_{n}}(w\{a_{3}(\lambda_{1})f^{2}-a_{2}(\lambda_{1})f^{1}\})+o_{P}(1).
\]

\end{document}